\documentclass{article}
\usepackage{mathrsfs}
\usepackage{ifthen}
\usepackage{amssymb}
\usepackage[all]{xy}
\usepackage[leftbars]{changebar}
\usepackage{vmargin}

\newcommand{\tiret}{\rule[0.6ex]{1.3ex}{0.26ex}}
\newcommand{\GL}[1]{\ensuremath{\mathrm{GL}}(#1)}
\newcommand{\Sscr}{\mathscr{S}}

\newcommand{\id}{\ensuremath{\mathrm{id}}}

\newcommand{\prods}[2]{{\langle #1,#2 \rangle}}
\newcommand{\ra}{\rightarrow}
\newcommand{\ZZ}{\mathbb{Z}}
\newcommand{\Cc}{\mathcal{C}}
\newcommand{\UU}{\mathbb{U}}
\newcommand{\CC}{\mathbb{C}}
\newcommand{\QQ}{\mathbb{Q}}
\newcommand{\al}{\alpha}
\newcommand{\Dscr}{\mathscr{D}}
\newcommand{\ga}{\gamma}
\newcommand{\ph}{\varphi}

\setpapersize{A4}
\setmarginsrb{20mm}{10mm}{20mm}{15mm}{10mm}{6mm}{0mm}{10mm}

\newcommand{\envfont}{\sf\bfseries}
\newcommand{\spacebeforeenv}{\vspace{1ex}}
\newcommand{\spaceafterenv}{\vspace{1ex}}

\newcounter{ctheo}
\renewcommand{\thectheo}{\arabic{ctheo}}

\newcommand{\metaenvironnementhm}[2]  
{
  \refstepcounter{ctheo}
  \ifthenelse{ \equal{#2}{toto} }{
    \spacebeforeenv \begin{traitsurlecote} \noindent {\envfont #1 \thectheo}
  }{
    \spacebeforeenv \begin{traitsurlecote} \noindent {\envfont #1 \thectheo{} \tiret\, #2.}
  }
}

\newenvironment{theo-eng}[1][toto]
{
  \metaenvironnementhm{Theorem}{#1}
}
{\end{traitsurlecote}\spaceafterenv}

\newenvironment{prop-eng}[1][toto]
{
  \metaenvironnementhm{Proposition}{#1}
}
{\end{traitsurlecote}\spaceafterenv}

\newenvironment{rem-eng}[1][toto]
{
  \metaenvironnementhm{Remark}{#1}
}
{\end{traitsurlecote}\spaceafterenv}

\newenvironment{lem-eng}[1][toto]
{ \metaenvironnementhm{Lemma}{#1} }
{\end{traitsurlecote}\spaceafterenv}

\newenvironment{cor-eng}[1][toto]
{ \metaenvironnementhm{Corollary}{#1} }
{\end{traitsurlecote}\spaceafterenv}

\newenvironment{dfn-eng}[1][toto]
{ \metaenvironnementhm{Definition}{#1} }
{\end{traitsurlecote}\spaceafterenv}

\newenvironment{notation}[1][toto]
{ \metaenvironnementhm{Notation}{#1} }
{\end{traitsurlecote}\spaceafterenv}

\newenvironment{demo-eng}
{\noindent  {\envfont Proof.~}}
{\spaceafterenv}

\newenvironment{traitsurlecote}
{\cbstart
\setcounter{changebargrey}{0}      
}
{\cbend
}

\newenvironment{rs}{\begin{list}{\tiret~}{
 \topsep=0.3ex
 \itemsep=0.3ex
 \labelsep=0em
 \parsep=0em
 \listparindent=1em
 \itemindent=0em
 \settowidth{\labelwidth}{--~}
 \leftmargin=\labelwidth
}}{\end{list}}

\newcommand{\numtoi}[1]{
  \ifthenelse{ \equal{#1}{1} }{i}{
  \ifthenelse{ \equal{#1}{2} }{ii}{
  \ifthenelse{ \equal{#1}{3} }{iii}{
  \ifthenelse{ \equal{#1}{4} }{iv}{
  \ifthenelse{ \equal{#1}{5} }{v}{
  \ifthenelse{ \equal{#1}{6} }{vi}{
  \ifthenelse{ \equal{#1}{7} }{vii}{
  \ifthenelse{ \equal{#1}{8} }{viii}{
  \ifthenelse{ \equal{#1}{9} }{ix}{
  \ifthenelse{ \equal{#1}{10} }{x}{
  \ifthenelse{ \equal{#1}{11} }{xi}{
  \ifthenelse{ \equal{#1}{12} }{xii}{
  \ifthenelse{ \equal{#1}{13} }{xiii}{
  \ifthenelse{ \equal{#1}{14} }{xiv}{
  \ifthenelse{ \equal{#1}{15} }{xv}{
  \ifthenelse{ \equal{#1}{16} }{xvi}{
  \ifthenelse{ \equal{#1}{17} }{xvii}{
  \ifthenelse{ \equal{#1}{18} }{xviii}{
    ERREUR,~MODIFIER~LA~MACRO~numtoi
  } } } } } } } } } } } } } } } } } }
}

\newcounter{cretraiti}
\newenvironment{ri}[1]{
\begin{list}{($\numtoi{\thecretraiti}$)}{
  \usecounter{cretraiti}
  \topsep=0.5ex
  \itemsep=0.3ex
  \labelsep=0.3em
  \parsep=0ex
  \listparindent=1em
  \settowidth{\labelwidth}{(#1)}
  \leftmargin=\labelwidth
}}{\end{list}
}

\newenvironment{keyword}
{\noindent  {\envfont Keywords.~}}
{\spaceafterenv}

\newenvironment{msc}
{\noindent  {\envfont Mathematics Subject Classification (2010).~}}
{\spaceafterenv}

\newcommand{\Oo}{\mathcal{O}}
\newcommand{\Hyp}{\mathscr{H}}
\newcommand{\Hscr}{\Hyp}
\newcommand{\Vreg}{\ensuremath{V^{\textrm{{\scriptsize reg}}}}}
\newcommand{\Pab}{\ensuremath{P^{\textrm{{\scriptsize ab}}}}}
\newcommand{\Ind}[3]{\ensuremath{\mathrm{Ind}_{#1}^{#2}(#3)}}
\newcommand{\Hom}[3]{\hom{#1}{#2}{#3}}
\renewcommand{\hom}[3]{\ensuremath{\mathrm{Hom}_{#1}(#2,#3)}}
\newcommand{\extensionzu}[5]{\ensuremath{\xymatrix{0 \ar[r]& #1 \ar^-{#4}[r] & #2 \ar^-{#5}[r] & #3\ar[r] & 1}}}

\newcommand{\extensionuu}[5]{\ensuremath{\xymatrix{1 \ar[r]& #1 \ar^-{#4}[r] & #2 \ar^-{#5}[r] & #3\ar[r] & 1}}}

\begin{document}

\title{Abelianization of Subgroups of Reflection Group and their Braid Group; an Application to Cohomology}

\author{Vincent Beck$^{1,2}$}

\maketitle

\footnotetext[1]{CMLA, ENS Cachan, CNRS, UniverSud, 61 Avenue du President Wilson, F-94230 Cachan}
\footnotetext[2]{IMJ, Universit\'e Paris 7, CNRS, 175 rue du Chevaleret, 75013 Paris}

\begin{abstract}
The final result of this article gives the order of the extension
$$\xymatrix{1\ar[r]  & P/[P,P] \ar^{j}[r] & B/[P,P] \ar^-{p}[r] & W \ar[r] & 1}$$
as an element of the cohomology group $H^2(W,P/[P,P])$ (where $B$ and $P$ stands for the braid group and the pure braid group associated to the complex
reflection group $W$). To obtain this result, we describe the abelianization of the stabilizer $N_H$ of a hyperplane $H$.
Contrary to the case of Coxeter groups, $N_H$ is not in general a reflection subgroup of the complex reflection group $W$.
So the first step is to refine Stanley-Springer's theorem on the abelianization of a reflection group.
The second step is to describe the abelianization of various types of big subgroups of the braid group $B$ of $W$.
More precisely, we just need a group homomorphism from the inverse image of $N_H$ by $p$ with values in $\QQ$ (where $p : B \ra W$ is the canonical morphism)
but a slight enhancement gives a complete description of the abelianization of $p^{-1}(W')$
where $W'$ is a reflection subgroup of $W$ or the stabilizer of a hyperplane. 
We also suggest a lifting construction for every element of the centralizer of a reflection in $W$.
\end{abstract}


\begin{msc} 20F36 - 20F55 - 20E22 \end{msc}

\section*{Introduction}\label{sec-intro}

Let us start with setting the framework. Let $V$ be a finite dimensional complex vector space;
a {\it reflection} is a non trivial finite order element $s$ of $\GL{V}$ which pointwise fixes a hyperplane of $V$, called {\it the hyperplane of $s$}.
The {\it line of $s$} is the one dimensional eigenspace of $s$ associated to the non trivial
eigenvalue of $s$. Let $W \subset \GL{V}$ be a {\it (complex) reflection group} that is to say a finite group generated by reflections.
We denote by $\Sscr$ the set of reflections of $W$ and $\Hyp$ the set of hyperplanes of $W$ :
$$\Sscr = \{s \in W,\quad \dim \ker(s- \id)= \dim V -1 \} \qquad \textrm{and } \qquad \Hyp = \{ \ker(s -\id),\quad s \in \Sscr \}\,.$$
For a reflection group, we denote by $\Vreg = V \setminus \cup_{H \in \Hyp} H$ the set of regular vectors.
According to a classical result of Steinberg~\cite[Corollary 1.6]{steinberg}, $\Vreg$ is precisely the set of vectors that no non trivial element of $W$ fixes.
Thus, the canonical map $\pi : \Vreg \ra \Vreg/W$
is a Galois covering. So let us fix a base point $x_0 \in \Vreg$ and denote by $P = \pi_1(\Vreg,x_0)$ and $B= \pi_1(\Vreg/W, \pi(x_0))$ the fundamental
groups of $\Vreg$ and its quotient $\Vreg/W$, we obtain the short exact sequence
\begin{equation}\label{PBW}
\xymatrix{1\ar[r]  & P \ar[r] & B \ar^{p}[r] & W \ar[r] & 1}
\end{equation}
\noindent The groups $B$ and $P$ are respectively called the {\it braid group} and the {\it pure braid group} of $W$.

The final result of this article (Corollary~\ref{cor-order-cocycle}) gives the order of the extension
$$\xymatrix{1\ar[r]  & P/[P,P] \ar^{j}[r] & B/[P,P] \ar^-{p}[r] & W \ar[r] & 1}$$
as an element of the cohomology group $H^2(W,P/[P,P])$. This order turns out to be the integer $\kappa(W)$ defined by Marin in~\cite{marin}. As explained in~\cite{marin},
this integer is linked to the periodicity of the monodromy representation of $B$ associated to the action of $W$ on its set of hyperplanes.

To obtain this result, we first describe in section~\ref{sec-reflection-group-case} the abelianization of some subgroups of complex reflection groups.
Specifically, we study the stabilizer of a hyperplane $H$ which is the same as the centralizer of a reflection of hyperplane~$H$.
Contrary to the case of Coxeter groups, this is not a reflection subgroup of the complex reflection group $W$ in general.
The first step is to refine Stanley-Springer's theorem~\cite{springer-lnm}\cite{stanley} on the abelianization
of a reflection group (see Proposition~\ref{prop-abel-stabilizer} in Section~\ref{sec-reflection-group-case}).

The rationale also relies on a good description (Section~\ref{sec-braid-group-case}) of abelianization of various types of big subgroups of the braid group $B$ of $W$
(here ``big'' stands for ``containing the pure braid group $P$'').
Though we just need to construct a group homomorphism from the inverse image of the stabilizer of $H$ by $p$ with values in $\QQ$ 
(see~Definition~\ref{dfn-morphism} and Proposition~\ref{prop-lin-caractere}),
we give in fact a complete description of the abelianization of $p^{-1}(W')$ with $W'$ a reflection subgroup of $W$ or the stabilizer of a hyperplane
(Proposition~\ref{prop-abelianisation} and Proposition~\ref{prop-abel-braid}). 
We also suggest a lifting construction for every element of the centralizer of a reflection in $W$
generalizing the construction of the generator of monodromy of~\cite[p.14]{bmr} (see Remark~\ref{rem-lifting}). 

Orbits of hyperplanes and ramification index are gathered in tables in the last section. 

We finish the introduction with some other notations.
We fix a $W$-invariant hermitian product on $V$ denoted by $\prods{\cdot}{\cdot}$ : orthogonality will always be relative to this particular hermitian product.
We denote by $S(V^*)$ the symmetric algebra of the dual of $V$ which is also the polynomial functions on $V$.

\begin{notation}[Around a Hyperplane of $W$]\label{notation-autour-hyperplan} For $H \in \Hyp$,
\begin{rs}
\item one chooses $\alpha_H \in V^*$ a linear form with kernel $H$;

\item one sets $W_H = \textrm{Fix}_{W}(H) = \{g \in W, \quad \forall\,x \in H, \quad gx = x \}$.
      This is a cyclic subgroup of $W$. We denote by $e_H$ its order and by $s_H$ its generator with determinant $\zeta_H = \exp(2i\pi/e_H)$.
      Except for identity, the elements of $W_H$ are precisely the reflections of $W$ whose hyperplane is $H$.
      The reflection $s_H$ is called the {\it distinguished reflection for $H$} in $W$.
\end{rs}
\end{notation}

For $n$ a positive integer, we denote by $\UU_n$ the group of the $n^{\textrm{{\scriptsize th}}}$ root of unity in $\CC$ and by $\UU$ the group of unit complex numbers. For a group $G$, we denote by $[G,G]$ the commutator subgroup of $G$ and by $G^{\textrm{{\scriptsize ab}}}= G/[G,G]$ the abelianization of $G$.

\vskip0.5ex
I wish to thank Ivan Marin for fruitful discussions and supports.

\section{Abelianization of Subgroups of Reflection Groups}\label{sec-reflection-group-case}

Stanley-Springer's Theorem~\cite{springer-lnm}\cite{stanley} gives an explicit description of the group of linear characters of
a reflection group using the conjugacy classes of hyperplanes. Naturally, it applies to all reflection subgroups of a reflection group
and in particular to parabolic subgroups thanks to Steinberg's theorem~\cite[Theorem 1.5]{steinberg}.
But since $P/[P,P]$ is the $\ZZ W$ permutation module defined by the hyperplanes of $W$,
we are interested in the stabilizer of a hyperplane which is not in general a reflection subgroup of $W$. So we have to go deeper in the study of the stabilizer of a
hyperplane.

For $H \in \Hscr$, we set $N_H = \{w \in W,\quad wH=H\} = \{w \in W,\quad  ws_H=s_Hw\}$ the stabilizer of $H$ which is also the centralizer of $s_H$.
We denote by $D = H^{\perp}$ the line of $s_H$ (or of every reflection of $W$ with hyperplane $H$) it is the unique $N_H$-stable line of $V$ such that $H \oplus D=V$
and we have $N_H=\{w \in W, \quad wD=D\}$ (see~\cite[Proposition 1.19]{broue}).
We denote by $C_H= \{w \in W,\quad  \forall\,x\in D, \ wx=x\}$ the parabolic subgroup associated to $D$.

Since $D$ is a line stable by every element of $N_H$ and $W_H \subset N_H$,
there exists an integer $f_H$ such that $e_H \mid f_H$ and the following sequence is exact
$$\extensionuu{C_H}{N_H}{\UU_{f_H}}{i}{r}$$
\noindent where $i$ is the natural inclusion and $r$  denote the restriction to $D$. We define $r$ to be the {\it natural linear character of $N_H$}.

Before stating our main result on the abelianization of the stabilizer of a hyperplane, let us start with a straightforward lemma.

\begin{lem-eng}[Abelianization of an exact sequence]\label{lem-abelianization} Let us consider the following exact sequence of groups
where $M$ is an abelian group.
$$\extensionuu{C}{N}{M}{i}{r}$$

Then the following sequence is exact
$$\xymatrix{C^\textrm{{\scriptsize ab}}\ar^{i^{\textrm{{\scriptsize ab}}}}[r]& N^{\textrm{{\scriptsize ab}}} \ar^{r^{\textrm{{\scriptsize ab}}}}[r]&
        M\ar[r]& 1}$$

Moreover the map $i^{\textrm{{\scriptsize ab}}}$ is injective if and only if $[N,N]= [C,C]$.
When $C^{\textrm{{\scriptsize ab}}}$ and $M$ are finite, the injectivity of $i^{\textrm{{\scriptsize ab}}}$ is also equivalent to $|C^{\textrm{{\scriptsize ab}}}||M| = |N^{\textrm{{\scriptsize ab}}}|$.

We also have the following criterion : the map $i^{\textrm{{\scriptsize ab}}}$ is injective if and only if the canonical restriction map
$\Hom{\textrm{{\scriptsize gr.}}}{N}{\CC^\times} \rightarrow \Hom{\textrm{{\scriptsize gr.}}}{C}{\CC^\times}$ is surjective.
\end{lem-eng}

\begin{demo-eng} We denote by $\theta$ the natural map from a group onto its abelianization. The map $r^{\textrm{{\scriptsize ab}}}$ is trivially surjective and
$r^{\textrm{{\scriptsize ab}}} i^{\textrm{{\scriptsize ab}}}$ is the trivial map. So let us consider $\theta(x) \in N^{\textrm{{\scriptsize ab}}}$ such that
$r^{\textrm{{\scriptsize ab}}}(\theta(x)) = 1$.
So we have $\theta(r(x)) = 1$. But $M$ is abelian and thus $r(x) = 0$. We then obtain that $x=i(y)$ for some $y \in C$.
Finally, we have $\theta(x) = i^{\textrm{{\scriptsize ab}}}(\theta(y))$.

The kernel of $i^{\textrm{{\scriptsize ab}}}$ is the image by the canonical map $\theta$ of $[N,N] \cap C$. But since $M$ is abelian, we have $[N,N] \subset C$.
So $i^{\textrm{{\scriptsize ab}}}$ is injective if and only if $[N,N] \subset [C,C]$ which is the same as $[N,N]=[C,C]$.

If we assume that $C^{\textrm{{\scriptsize ab}}}$ and $M$ are finite
then $|N^{\textrm{{\scriptsize ab}}}| = |i^{\textrm{{\scriptsize ab}}} (C^{\textrm{{\scriptsize ab}}})||M|$.
So $|C^{\textrm{{\scriptsize ab}}}||M| = |N^{\textrm{{\scriptsize ab}}}|$ is equivalent to $|i^{\textrm{{\scriptsize ab}}}(C^{\textrm{{\scriptsize ab}}})| = |C^{\textrm{{\scriptsize ab}}}|$ which is equivalent to $i^{\textrm{{\scriptsize ab}}}$ is injective.

Since $\CC^\times$ is a commutative group, we have the following commutative square whose vertical arrows are isomorphisms
$$\xymatrix{\Hom{\textrm{{\scriptsize gr.}}}{N}{\CC^\times}\ar[r] \ar[d] & \Hom{\textrm{{\scriptsize gr.}}}{C}{\CC^\times} \ar[d]\\
\Hom{\textrm{{\scriptsize gr.}}}{N^{\textrm{{\scriptsize ab}}}}{\CC^\times} \ar^{\circ i^{\textrm{{\scriptsize ab}}}}[r]&
\Hom{\textrm{{\scriptsize gr.}}}{C^{\textrm{{\scriptsize ab}}}}{\CC^\times}}$$
Moreover, since $\CC^\times$ is a divisible abelian group, $i^{\textrm{{\scriptsize ab}}}$ is injective if and only if $\circ i^{\textrm{{\scriptsize ab}}}$
is surjective.
\end{demo-eng}

Before applying the preceding lemma to the stabilizer of a hyperplane, let us introduce a definition

\begin{dfn-eng}[Commuting hyperplane]\label{dfn-commuting-hyperplane} Let $H,H' \in \Hscr$. We say that {\it $H$ and $H'$ commute} if
$s_H$ and $s_{H'}$ commute.

We have the following equivalent characterizations (see~\cite[Lemma 1.7]{broue}):
\begin{ri}{(iii)}
\item the hyperplanes $H$ and $H'$ commute
\item $H=H'$ or $D = H^\perp \subset H'$
\item every reflection of $W$ with hyperplane $H$ commutes with every reflection of $W$ with hyperplane $H'$
\item there exists a reflection of $W$ with hyperplane $H$ which commutes with a reflection of $W$ with hyperplane~$H'$
\end{ri}

We denote by $\Hscr_H$ the set of hyperplanes commuting with $H$ and $\Hscr'_H = \Hscr_H \setminus \{ H \}$.
\end{dfn-eng}

As a consequence of Lemma~\ref{lem-abelianization}, we are now able to formulate the following proposition.

\begin{prop-eng}[Abelianization of the stabilizer of a hyperplane]\label{prop-abel-stabilizer} For a hyperplane $H \in \Hscr$, the following sequence is exact
$$\xymatrix{C_H^\textrm{{\scriptsize ab}}\ar^{i^{\textrm{{\scriptsize ab}}}}[r]& N_H^{\textrm{{\scriptsize ab}}} \ar^{r^{\textrm{{\scriptsize ab}}}}[r]&
        \UU_{f_H}\ar[r]& 1}$$

\noindent Moreover, we have the following geometric characterization of the injectivity of $i^{\textrm{{\scriptsize ab}}}$ :
the map $i^{\textrm{{\scriptsize ab}}}$ is injective if and only if the orbits of the hyperplanes commuting with $H$ under $N_H$ and $C_H$ are the same.

This condition is verified for every hyperplane of every complex reflection group except the hyperplanes of the exceptional group $G_{25}$ and
the hyperplanes $H_i$ ($1 \leq i \leq r$) of the group $G(de,e,r)$ when $r=3$ and $e$ is even (see section~\ref{sec-table} for the notations).
\end{prop-eng}

\begin{demo-eng} Let us assume that the orbits of the hyperplanes commuting with $H$ under $N_H$ and $C_H$ are the same.
Steinberg's theorem and definition~\ref{dfn-commuting-hyperplane}~$(ii)$ ensure us that $C_H$ is the reflection subgroup of $W$ generated
by the $s_{H'}$ for $H' \neq H$ commuting with $H$.
Thanks to Stanley-Springer's theorem, we are able to describe the linear characters of $C_H$. For every linear character $\delta$ of $C_H$,
there exists integers $e_\Oo$ for $\Oo \in \Hscr_H'/C_H$ such that the polynomial
$$Q_\delta = \prod_{\Oo \in \Hscr_H'/C_H} \prod_{H' \in \Oo} {\al_{H'}}^{e_\Oo} \in S(V^*)$$
verifies $g Q_\delta = \delta(g)Q_\delta$ for every $g \in C_H$. Since the orbits of $\Hscr_H'$ under $N_H$ are the same as the orbits of
$\Hscr'_H$ under $C_H$, then for every $g \in N_H$ and every $\Oo \in \Hscr'_H/C_H$ there exists $\lambda_g \in \CC^\times$ such that
$$g \prod_{H' \in \Oo} {\al_{H'}}= \lambda_g  \prod_{H' \in \Oo} {\al_{H'}}$$
So we obtain that, for every $g \in N_H$, there exists $\mu_g \in \CC^\times$ such that $g Q_\delta = \mu_g Q_\delta$.
Thus every linear character of $C_H$ extends to $N_H$ and Lemma~\ref{lem-abelianization} tells us that the map $i^{\textrm{{\scriptsize ab}}}$ is injective.

Let us assume now that every linear character of $C_H$ extends to $N_H$. The orbit of $H$ under $N_H$ and $C_H$ is $\{ H \}$.
So let us consider an orbit $\Oo \in \Hscr'_H/C_H$. We define
$$Q = \prod_{H' \in \Oo} \al_{H'} \in S(V^*)\,.$$
Thanks to Stanley-Springer's theorem, $Q$ define a linear character $\chi$ of $C_H$ : for every $c \in C_H$, there exists $\chi(c) \in \CC^\times$ such that
$cQ = \chi(c)Q$ for every $c \in C_H$. We then consider the $N_H$-submodule $M$ of $S(V^*)$ generated by $Q$.
As a vector space, $M$ is generated by the family $(nQ)_{n \in N_H}$. But, since $C_H$ is normal in $N_H$, we obtain for $c \in C_H$,
$$cnQ = nn^{-1}cn Q= n \chi(n^{-1}cn)Q\,.$$
Since $\chi$ extends to $N_H$, we have $\chi(n^{-1}cn) = \chi(c)$ and then $cnQ=\chi(c)nQ$. Stanley-Springer's theorem allows us to conclude that
$nQ= \lambda_n Q$ for some $\lambda_n \in \CC^\times$. Since $S(V^*)$ is a UFD, we obtain that $\Oo$ is still an orbit under $N_H$.

In Section~\ref{sec-table}, we give tables for the various orbits of hyperplanes for the infinite series $G(de,e,r)$.
For the exceptional complex reflection groups, we check the injectivity of $i^{\textrm{{\scriptsize ab}}}$ using the package {\sf CHEVIE} of {\sf GAP}
\cite{chevie}\cite{gap} : see~Remark~\ref{rem-gap} in Section~\ref{sec-table}.
\end{demo-eng}

For a hyperplane $H \in \Hscr$, the comparison of $e_H$ and $f_H$ leads to the following definition.

\begin{dfn-eng}[Ramification at a hyperplane] We define $d_{H} = f_{H}/e_{H}$ to be {\it the index of ramification of $W$ at the hyperplane $H$}.
We say that $W$ is {\it unramified at $H$} if $d_{H}=1$.

We say that an element $w \in N_H$ such that $r(w) = \exp(2i\pi/f_{H})$ {\it realizes the ramification}.
\end{dfn-eng}

\begin{rem-eng}[The Coxeter Case]
When $H$ is unramified, we have $N_H = C_H \times W_H$ which is generated by reflections thanks to Steinberg's theorem~\cite[Theorem 1.5]{steinberg} and
$s_{H}$ realizes the ramification. Moreover $i^{\textrm{\scriptsize ab}}$ is trivially injective.

In a Coxeter group, every hyperplane is unramified. Indeed, the eigenvalue on the line $D$ of an element of $N_H$ is a finite order element of the field of the real numbers.
\end{rem-eng}

In section~\ref{sec-table}, we give tables for the values of $e_H,f_H$ and $d_H$ for every hyperplane of every complex reflection groups.
From these tables, we obtain the following remarks.

\begin{rem-eng}[Unramified $G(de,e,r)$]
All the hyperplanes of $G(de,e,r)$ are unramified only when $r=1$ or when $G(de,e,r)$ is a Coxeter group (that is to say if $d=2$ and $e=1$ and $r \geq 2$
(Coxeter group of type $B_r$) or if $d=1$ and $e=2$ and $r \geq 3$ (Coxeter group of type $D_r$) or if $d=1$ and $e=1$ and $r\geq 3$ (Coxeter group
of type $A_{r-1}$) of if $d=1$ and $r=2$ (Coxeter group of type $I_2(e)$).
\end{rem-eng}

\begin{rem-eng}[Unramified exceptional groups] The only non Coxeter groups for which every hyperplane is unramified
are $G_8,G_{12}$ and $G_{24}$.
\end{rem-eng}

\section{Abelianization of Subgroups of Braid Groups}\label{sec-braid-group-case}

In this section, we describe abelianizations of subgroups of $B$ containing $P$ that is to say of inverse images of subgroups $W'$ of $W$.
Explicitly, we are able to give a complete description of ${p^{-1}(W')}^{\textrm{{\scriptsize ab}}}$ if $W'$ is a reflection subgroup (Proposition~\ref{prop-abelianisation})
or if $W'$ is the stabilizer of a hyperplane under geometrical assumptions on the hyperplane (Proposition~\ref{prop-abel-braid}). 
We also construct a particular linear character of $p^{-1}(N_H)$ lifting 
the natural linear character $r$ of $N_H$ which is of importance for the next section (Definition~\ref{dfn-morphism}).

Our method is similar to the method of~\cite{bmr} for the description of $B^{\textrm{\scriptsize ab}}$ :
we integrate along pathes invariants polynomial functions. So, we have first to construct invariant polynomial functions and then
verify that we have constructed enough of them.

\subsection{Subgroup Generated by Reflections}

In this subsection, we fix $C$ a subgroup of $W$ generated by reflections.
We denote by $\Hyp_C\subset \Hscr$ the set of hyperplanes of $C$. For $H \in \Hyp_C$, then $C' = \{c \in C, \quad \forall\,x \in H, \ cx=x\}$ is a subgroup of $W_H$
and so generated by ${s_H}^{a_H}$ with $a_H \mid e_H$. For $H \in \Hyp \setminus \Hyp_C$, we set $a_H= e_H$.
Then $C= \langle {s_H}^{a_H}, \ H \in \Hyp\rangle$. For $\Cc \in \Hyp/C$ a $C$-class of hyperplanes of $W$,
we denote by $a_\Cc$ the common value of $a_H$ for $H \in \Cc$.

The aim of this subsection is to give a description of the abelianization of $p^{-1}(C) \subset B$. For this, we follow the method of~\cite{bmr}
and we start to exhibit invariants which will be useful to show the freeness of our generating set of $p^{-1}(C)^{\textrm{{\scriptsize ab}}}$.

\begin{lem-eng}[An invariant]\label{lem-invariant}
We define, for $\Cc\in \Hyp/C$,
$$\alpha_\Cc = \prod_{H \in \Cc} {\alpha_H}^{e_H/a_H} \in S(V^*)\,.$$
Then $\alpha_\Cc$ is invariant under the action of $C$.
\end{lem-eng}

\begin{demo-eng} If $\Cc$ is a class of hyperplanes of $\Hyp_C$ then this is an easy consequence of Stanley-Springer's theorem.

Assume that $\Cc$ is not a class of hyperplanes of $\Hyp_C$. Let us choose a reflection $s$ of $C$ and let $n_s$ be the order of $s$.
Since $\Cc$ is not a class in $\Hyp_C$, the hyperplane of $s$ does not belong to $\Cc$.
We then deduce that the orbits of $\Cc$ under the action of $\langle s \rangle$ are of two types.

First type : $H \in \Cc$ such that $s_Hs = ss_H$. We denote by $H_s$ and $D_s$ the hyperplane and the line of $s$.
Since $H\neq H_s$, Lemma 1.7 of~\cite{broue} tells us that
$D=H^\perp \subset H_s$ and $D_s \subset H$. For $x \in V$, we write $x=  h+d$ with $h \in H$ and $d \in D$.
Since $D_s \subset H$, $H$ is stable by $s$ and we have $s^{-1}(x) = s^{-1}(h+d) = h' + d$ with $h' \in H$. So we obtain
$s\alpha_H(x) = \alpha_H(d) = \alpha_H(x)$ and the orbit of $H$ under $\langle s\rangle$ is reduced to $H$ and moreover $s\alpha_H = \alpha_H$.

Second type : $H \in \Cc$ such that $s_H s \neq s s_H$. If $s^iH=H$ then $s^i$ and $s_H$ commute and thus, Lemma 1.7 of~\cite{broue} ensures us that $s^{i}$ is trivial.
We then obtain that the orbit of $H$ under $\langle s \rangle$ has cardinal $n_s$. So if we write
$Q:=\alpha_H \alpha_{sH}\cdots \alpha_{s^{n_s-1}H} = \lambda \alpha_H s\alpha_{H} \cdots s^{n_s-1}\alpha_{H}$ with $\lambda \in \CC^{\times}$, we have $sQ = Q$.

We then easily obtain $s\alpha_\Cc=\alpha_\Cc$ for every $s \in C$ and so $\alpha_\Cc$ is invariant under the action of $C$.
\end{demo-eng}

Before stating the main result of the subsection, we recall the notion of ``a generator of the monodromy around a hyperplane'' as explained in~\cite[p.14]{bmr}.
For $H \in \Hyp$, we define a generator of the monodromy around $H$ to be a path $s_{H,\ga}$ in $\Vreg$ which is the composition of three pathes.
The first path is a path $\ga$ going from $x_0$ to a point $x_H$ which is near $H$ and far from other hyperplanes.
To describe the second path, we write $x_H = h+d$ with $h \in H$ and $d \in D=H^\perp$, and the second path is
$t \in [0,1] \mapsto h +\exp(2 i \pi t/e_H)d$ going from $x_H$ to $s_H(x_H)$. The third path is $s_H(\ga^{-1})$ going from $s_H(x_H)$ to $s_H(x_0)$.
We can now state our abelianization result.

\begin{prop-eng}[Abelianization of subgroups of the braid group]\label{prop-abelianisation} Let $C$ be a subgroup of $W$ generated by reflections.
Then $p^{-1}(C)^{\textrm{\scriptsize ab}}$ is the free abelian group over $\Hyp/C$ the $C$-classes of hyperplanes of $W$.

More precisely, we have $p^{-1}(C) = \langle {s_{H,\ga}}^{a_H},\ (H,\gamma)  \rangle$ (see~\cite[Theorem 2.18]{bmr}).
For $\Cc \in \Hyp/C$, we denote by $(s_\Cc^{a_\Cc})^{\textrm{{\scriptsize ab}}}$ the common value in $p^{-1}(C)^{\textrm{{\scriptsize ab}}}$
of the ${s_{H,\ga}}^{a_H}$ for $H \in \Cc$.
Then $p^{-1}(C)^{\textrm{{\scriptsize ab}}} = \langle (s_\Cc^{a_\Cc})^{\textrm{{\scriptsize ab}}},\quad \Cc \in \Hyp/C \rangle$.
Moreover, for $\Cc \in \Hyp/C$, there exists a group homomorphism $\ph_\Cc : p^{-1}(C) \ra \ZZ$ such that $\ph_\Cc((s_\Cc^{a_\Cc})^{\textrm{{\scriptsize ab}}})=1$ and
$\ph_\Cc((s_{\Cc'}^{a_{\Cc'}})^{\textrm{{\scriptsize ab}}})=0$ for $\Cc' \neq \Cc$.
\end{prop-eng}

\begin{demo-eng} First of all, Lemma~2.14.(2) of~\cite{bmr} shows that ${s_{H,\ga}}^{a_H} = {s_{H,\ga'}}^{a_H}$ in $p^{-1}(C)^{\textrm{{\scriptsize ab}}}$.
Now, for $c \in C$, we choose $x \in p^{-1}(C)$ such that $p(x) = c$. We have $xs_{H,\ga}x^{-1}= s_{cH,x (c\ga)}$.
So ${s_{cH,x (c\ga)}}^{a_{cH}}$ and ${s_{H,\ga}}^{a_H}$ are conjugate by an element of $p^{-1}(C)$. So, we have
$${s_{cH,x (c\ga)}}^{a_{cH}} = {s_{H,\ga}}^{a_H} \in p^{-1}(C)^{\textrm{{\scriptsize ab}}}\,.$$
\noindent And then $p^{-1}(C)^{\textrm{{\scriptsize ab}}} = \langle (s_\Cc^{a_\Cc})^{\textrm{{\scriptsize ab}}},\quad \Cc \in \Hyp/C \rangle$.

Let us now show that the family $((s_\Cc^{a_\Cc})^{\textrm{{\scriptsize ab}}})_{\Cc \in \Hyp/C}$ is free over $\ZZ$.
We identify $p^{-1}(C)$ with
$$p^{-1}(C) = \left(\displaystyle \bigsqcup_{c,c' \in C} \pi_1(\Vreg, cx_0,c'x_0) \right)/ C$$
\noindent where $\pi_1(\Vreg, cx_0,c'x_0)$ denotes the homotopy classes of pathes from $c(x_0)$ to $c'(x_0)$ and the action of $C$ on path is simply the composition.

Since $\alpha_\Cc : \Vreg \ra \CC^{\times}$ is $C$-invariant (Lemma~\ref{lem-invariant}), the functoriality of $\pi_1$ defines a group homomorphism
$\pi_1(\alpha_\Cc)$ from $p^{-1}(C)$ to $\pi_1(\CC^\times, \alpha_\Cc(x_0))$. Moreover, the map
$$I: \ga \longmapsto \frac{1}{2i\pi}\int_{\ga} \frac{\textrm{d} z}{z}$$
\noindent realizes a group isomorphism between $\pi_1(\CC^\times, \alpha_\Cc(x_0))$ and $\ZZ$. The composition of this two maps defines the group homomorphism $\ph_\Cc$.

For $H \in \Cc$ and $\Cc' \in \Hyp/C$, let us compute $\ph_{\Cc'}({s_{H,\ga}}^{a_H})$.
The path ${s_{H,\ga}}^{a_H}$ is the composition of three paths. The first one is $\ga$, the third one is ${s_H}^{a_H}(\ga^{-1})$ and the second one is
$\eta : t \in [0,1] \mapsto h +\exp(2 i \pi a_H t/e_H)d$.

Since $\alpha_{\Cc'}$ is $C$-invariant, when we apply $\pi_1(\alpha_{\Cc'})$, the first part of the path and the third one are inverse from each other.
So when applying $I$, they do not appear. We thus obtain
$$\ph_{\Cc'}({s_{H,\ga}}^{a_H})= \frac{1}{2i\pi} \int_{\alpha_{\Cc'} \circ \eta} \frac{\textrm{d} z}{z}\,.$$
\noindent Using the logarithmic derivative, we obtain
$$\ph_{\Cc'}({s_{H,\ga}}^{a_H})= \frac{1}{2i\pi} \sum_{H'\in \Cc'} \frac{e_{H'}}{a_{H'}}
    \int_{0}^{1} 2i\pi \frac{a_H}{e_H} \exp(2i\pi a_H t/e_H) \frac{\alpha_{H'}(d)}{\alpha_{H'}(h+\exp(2i\pi a_Ht/e_H)d)} \textrm{d} t$$
\noindent To compute this sum, we regroup the terms after the orbit of $H'$ under $\langle {s_H}^{a_H} \rangle$.

Lemma~1.7 of~\cite{broue} shows that there are three types of orbits : two types of orbits reduced to one single hyperplane and
one other type of orbits corresponding to reflections that do not commute with $s_H$.

Let us first study the orbits reduced to one single hyperplane. The first type corresponds to the hyperplane $H$ whose term of the sum is $1$ and
this term appears if and only if $H \in \Cc'$. The second type corresponds to hyperplanes $H'$ such that $D = H^{\perp} \subset H'$.
The corresponding term of the sum is $0$ since $\alpha_{H'}(d) =0$

Let us now study the other orbits.
The orbits of $H'$ under ${s_H}^{a_H}$ is $\{ H', \ldots, {{s_H}^{a_H}}^{(e_H/a_H-1)}(H')\}$. Moreover,
a quotient of the form $\alpha_{H'}(x)/\alpha_{H'}(y)$ does not depend of the linear form with kernel $H'$. So we can replace $\alpha_{{s_H}^iH'}$ by $s_H^i\alpha_{H'}$ to obtain
$$\begin{array}{rcl} \displaystyle \frac{\exp(2i\pi a_H t/e_H)\alpha_{{s_H}^{-a_Hk}H'}(d)}{\alpha_{{s_H}^{-a_Hk}H'}(h) + \exp(2i\pi a_H t/e_H)\alpha_{{s_H}^{-a_Hk}H'}(d)} &=&
\displaystyle \frac{\exp(2i\pi a_H t/e_H){s_H}^{-a_Hk}\alpha_{H'}(d)}{{s_H}^{-a_Hk}\alpha_{H'}(h) + \exp(2i\pi a_H t/e_H){s_H}^{-a_Hk}\alpha_{H'}(d)}\\[3ex] &=& \displaystyle
        \frac{\exp(2 i \pi a_H (t+k)/e_H)\alpha_{H'}(d)}{\alpha_{H'}(h) + \exp(2 i \pi a_H (t+k)/e_H)\alpha_{H'}(d)}\,. \end{array}$$
Considering the sum over the orbit under $\langle {s_H}^{a_H} \rangle$ of $H'$, we obtain
$$\sum_{k=0}^{e_H/a_H-1} \int_{0}^{1} \frac{\exp(2 i \pi a_H (t+k)/e_H)\alpha_{H'}(d)}{\alpha_{H'}(h) + \exp(2 i \pi a_H (t+k)/e_H)\alpha_{H'}(d)}\,\textrm{d} t =
        \frac{e_H}{a_H}\int_{0}^{1} \frac{\exp(2 i \pi t)\alpha_{H'}(d)}{\alpha_{H'}(h) + \exp(2 i \pi t)\alpha_{H'}(d)}\,\textrm{d} t\,.$$
\noindent Since $x_H$ is chosen such that $\alpha_{H'}(h) \neq 0$ for $H' \neq H$ and $d$ is small, the last term is $0$ as the index of the circle of center $0$
and radius $|\alpha_{H'}(d)|$ relatively to the point $-\alpha_{H'}(h)$.
\end{demo-eng}

\begin{rem-eng}[Extreme cases] The two extreme cases where $C = 1$ and $C=W$ may be found in~\cite[Prop. 2.2.(2)]{bmr} and~\cite[Theorem 2.17.(2)]{bmr}.
In the first case, $p^{-1}(C) = P$ is the pure braid group whose abelianization is the free abelian group over $\Hyp$.
In the second case $p^{-1}(C) = B$ is the braid group whose abelianization is the free abelian group over $\Hyp/W$.
\end{rem-eng}

\begin{rem-eng}\label{rem-prolongement} The logarithmic derivative shows that for every $\gamma \in p^{-1}(C)$ and $n \in \ZZ$, we have
$$\int_{{\al_\Cc}^{n}\circ \gamma} \frac{\textrm{d} z}{z} = n\ph_\Cc(\ga)\,.$$
\end{rem-eng}

\subsection{Stabilizer of a hyperplane}

Let us recall the notation of Section~\ref{sec-reflection-group-case}; we consider $H \in \Hscr$ a hyperplane of the reflection group $W$.
We denote by $N_H$ the stabilizer of $H$ in $W$ and $C_H$ the parabolic subgroup of $W$ associated to the line $D= H^\perp$.
The set of hyperplanes commuting with $H$ is $\Hscr_H$ (see definition~\ref{dfn-commuting-hyperplane}).

\subsubsection{A group homomorphism}

The aim of this paragraph is to construct an ``extension'' of the natural character of $N_H$ to the group $p^{-1}(N_H)$ which will be useful for the
third section. We still follow the method of~\cite{bmr} : we construct an invariant function with values in $\CC^\times$ (Lemma~\ref{lem-invariant2})
and integrate it (Definition~\ref{dfn-morphism}). To obtain the ``extension'' properties of the linear character of $p^{-1}(N_H)$ (Proposition~\ref{prop-lin-caractere}), we construct a lifting in the braid group of the elements of $N_H$ (Remark~\ref{rem-lifting}). 
This lifting is traced from the construction of the generator of the monodromy.

\begin{lem-eng}[The invariant function]\label{lem-invariant2} The element $\alpha_{N_H} = {\alpha_{H}}^{f_{H}} \in S(V^*)$ is invariant under $N_H$.
\end{lem-eng}

\begin{demo-eng} For $n \in N_H$ and $x = h+d \in V$ with $h \in H$ and $d \in D$,  we have
$$n\alpha_{H}(x) = \alpha_{H}(n^{-1}h + n^{-1}d) = \alpha_{H}(h'+ r(n^{-1})d)= r(n^{-1}) \alpha_{H}(d) = r(n^{-1})\alpha_{H}(x)$$
with $h' \in H$.
\end{demo-eng}

\begin{dfn-eng}[The group homomorphism]\label{dfn-morphism} As in the proof of Proposition~\ref{prop-abelianisation}, we write
$$p^{-1}(N_H) = \left(\displaystyle \bigsqcup_{n,n' \in N_H} \pi_1(\Vreg, nx_0,n'x_0) \right)/ N_H\,.$$
Since $\alpha_{N_H} : \Vreg \ra \CC^{\times}$ is $N_H$-invariant (Lemma~\ref{lem-invariant2}), the functoriality of $\pi_1$ allows us to define a group homomorphism
$\pi_1(\alpha_{N_H})$ from $p^{-1}(N_H)$ to $\pi_1(\CC^\times, \alpha_{N_H}(x_0))$. Moreover, the map
$$I: \ga \longmapsto \frac{1}{2i\pi}\int_{\ga} \frac{\textrm{d} z}{z}$$
\noindent realizes a group isomorphism between $\pi_1(\CC^\times, \alpha_{N_H}(x_0))$ and $\ZZ$. The composition of this two maps defines a group homomorphism
$\rho' : p^{-1}(N_H) \ra \ZZ$. We also define $\rho = {f_H}^{-1} \rho' : p^{-1}(N_H) \ra \QQ$.
\end{dfn-eng}

\begin{rem-eng}[Center of the braid group of $G_{31}$] In~\cite[Theorem 2.24]{bmr}, it is shown that the center of the braid group $B$ 
of an irreducible reflection group $W$ is an infinite cyclic group generated by $\beta : t \mapsto \exp(2i\pi t/|Z(W)|)x_0$
(where $x_0 \in \Vreg$ is a base point) for all but six exceptional reflection groups. In his articles~\cite{bessis1}\cite{bessis2},
Bessis proves that the result holds for all reflection groups but the exceptional one $G_{31}$.

This remark is a first step toward the case of $G_{31}$ : we show that if $ZB$ is an infinite cyclic group, it is generated by $\beta$. 
For this, let us consider $H \in \Hscr$ a hyperplane of~$G_{31}$ and $\rho'$ the group homomorphism defined above. 
Since $ZB \subset p^{-1}(N_H)$, $\rho'$ restricts to a group homomorphism from $ZB$ to $\ZZ$ such that $\rho'(\beta) = 1$.
So if $ZB$ is an infinite cyclic group, it is generated by $\beta$. 
\end{rem-eng}

\begin{rem-eng}[The lifting construction]\label{rem-lifting}
Let us consider $w \in N_H$. We now construct a path $\widetilde{w}$ in $\Vreg$ starting from $x_0$ and ending at $w(x_0)$ :
$p(\widetilde{w})= w$. We use the notations of the description of the generators of the monodromy around $H$ :
we write $x_H = h + d$ with $h \in H$ and $\al_{H'}(h) \neq 0$ for $H' \neq H$ and $d \in D = H^\perp$.
Since $w \in N_H$, we have $w(x_H) = h' + \exp(2i k\pi/f_H)d$ with $h' \in H$ and $0 \leq k < f_H$.

The path $\widetilde{w}$ consists into four parts. As in the case of the generators of the monodromy,
the first part is a path $\ga$ from $x_0$ to $x_H$ and the fourth path is $w(\ga^{-1})$ from $w(x_H)$ to $w(x_0)$.
Let us now describe the second part and the third part.
The second part of $\widetilde{w}$ is the path
$$t \in [0,1] \mapsto h + \exp(2ik\pi t/f_H)d \in \Vreg\,.$$
The third part is of the form $t \in [0,1] \mapsto \theta(t) + \exp(2ik\pi/f_H) d$ where $\theta(t)$ is a path in the complex affine line $\Dscr$ generated by
$h'$ and $h$. It is easy to force the third part of $\widetilde{w}$ to stay in $\Vreg$ since its image is contained in the affine line
$\exp(2ik\pi/f_H) d + \Dscr$ which is parallel to the hyperplane $H$ and meets each of the other hyperplanes in a single point :
so we just have to avoid a finite number of points in $\CC$.
\end{rem-eng}

It remains to show that the constructed group homomorphism $\rho$ is an ``extension'' of the natural character of $N_H$.
More precisely, we have the following proposition.

\begin{prop-eng}[The ``extension'' property]\label{prop-lin-caractere} We have the following commutative square
$$\xymatrix{p^{-1}(N_H)\ar^-{\rho}[r] \ar_{p}[d] & \QQ \ar^{\pi'}[d] \\ N_H \ar^-{r}[r]& \UU_{f_H}}$$
where $\pi' : x \in \QQ \mapsto \exp(2i\pi x)$.
\end{prop-eng}

\begin{demo-eng} Let us choose $w \in N_H$ realizing the ramification of $H$ that is to say such that $r(w)= \exp(2i\pi/f_H)$.
Steinberg's theorem~\cite[Theorem 1.5]{steinberg} tells us $N_H = \langle w, s_{H'},\ H' \in \Hyp_H \rangle$.
We consider $\widetilde{w} \in B$ such that $p(\widetilde{w}) =w$ constructed in Remark~\ref{rem-lifting}. Theorem~2.18 of~\cite{bmr} shows that
$$p^{-1}(N_H)= \langle \widetilde{w}, s_{H',\ga}, {s_{H'',\ga'}}^{e_{H''}},\quad H' \in \Hyp_H,\ H'' \in \Hyp \setminus \Hyp_H,\ \ga,\ \ga' \rangle\,.$$
As a consequence, we only have to verify the following equalities to obtain the commutative square
\begin{ri}{(iii)}
\item $\rho'(\widetilde{w})=1$
\item $\rho'(s_{{H},\ga}) = f_H/e_H$
\item $\rho'(s_{H',\ga}) = 0$ for $H' \in \Hyp_H' = \Hyp_H\setminus \{H \}$
\item $\rho'({s_{H',\ga}}^{e_{H'}}) = 0$ for $H' \in \Hyp \setminus \Hyp_H$.
\end{ri}

As in the proof of Proposition~\ref{prop-abelianisation}, the $\ga$-part of $s_{H,\ga}$
(resp. $s_{H',\ga}$ for $H'\in \Hyp_H'$ and ${s_{H',\ga}}^{e_{H'}}$ for $H' \in \Hyp \setminus \Hyp_H$) does not step in the computation of $\rho'$.
We thus obtain
$$\rho'(s_{H,\ga}) = \frac{1}{2i\pi}\int_{0}^{1} f_{H} \frac{2i\pi}{e_{H}}
        \frac{\alpha_{H}(\exp(2i\pi t/e_{H})d)}{\alpha_{H}(h+d\exp(2i\pi t/e_{H}))}\textrm{d}t = f_H/e_H\,.$$
For $H' \in \Hyp_H'$, we set $x_{H'} = h'+d'$ with $h' \in H'$ and $d' \in D'={H'}^\perp$. We then obtain
$$\rho'(s_{H',\ga}) = \frac{1}{2i\pi}\int_{0}^{1} f_{H} \frac{2i\pi}{e_{H'}}
        \frac{\alpha_{H}(\exp(2i\pi t/e_{H'})d')}{\alpha_{H}(h'+d'\exp(2i\pi t/e_{H'}))}\textrm{d}t = 0$$
\noindent since $\alpha_{H}(d')=0$ for $H \in \Hyp_H'$. With the same arguments, we obtain for $H' \in \Hyp \setminus \Hyp_H$
$$\rho'({s_{H',\ga}}^{e_H}) = \frac{1}{2i\pi}\int_{0}^{1} 2i\pi f_{H}
        \frac{\exp(2i\pi t)\alpha_{H}(d')}{\alpha_{H}(h')+\exp(2i\pi t)\alpha_{H}(d')}\textrm{d}t = 0$$
\noindent since $d'$ is small and $\alpha_{H}(h')\neq 0$.

For $\widetilde{w}$, neither the first and fourth part step in the computation nor the third one. Moreover,
as in the computation of $\rho'(s_{H,\ga})$ the second part of $\widetilde{w}$ gives $1$.
\end{demo-eng}

\subsubsection{The Stabilizer Case}

In this paragraph, we extend the results of Section~\ref{sec-reflection-group-case} to the braid group.
Namely, we have the following short exact sequence
$$\xymatrix{1 \ar[r] & p^{-1}(C_H) \ar^j[r] & p^{-1}(N_H) \ar^-{rp}[r] & \UU_{f_H} \ar[r]& 1}$$
which gives rise to the following exact sequence (Lemma~\ref{lem-abelianization})
$$\xymatrix{p^{-1}(C_H)^{\textrm{\scriptsize ab}} \ar^{j^{\textrm{\scriptsize ab}}}[r]
    & p^{-1}(N_H)^{\textrm{\scriptsize ab}} \ar^-{(rp)^{{\textrm{\scriptsize ab}}}}[r] & \UU_{f_H} \ar[r]& 1}$$
and Proposition~\ref{prop-abel-stabilizer} extends to the braid group in the following way.

\begin{prop-eng}[Abelianization in the braid group]\label{prop-abel-braid}
If the orbits of the hyperplanes of $\Hscr$ under $N_H$ and $C_H$ are the same, the map $j^{\textrm{{\scriptsize ab}}}$ is injective.

Moreover under this hypothesis, $p^{-1}(N_H)^{\textrm{\scriptsize ab}}$ is the free abelian group with basis
$\widetilde{w}$, $(s_\Cc)^{\textrm{\scriptsize ab}}$ for $\Cc \in \Hscr'_H/C_H$ and $(s_{\Cc}^{e_\Cc})^{\textrm{\scriptsize ab}}$
for $\Cc \in (\Hscr\setminus \Hscr_H)/C_H$.

Let us consider the infinite series (see Section~\ref{sec-table} for notations).
If $H=H_i$, the  orbits under $N_H$ and $C_H$ are always the same except when $r=3$ and $e$ is even and when $r=2$ and $e \geq 3$.
If $H=H_{i,j,\zeta}$, the orbits under $N_H$ and $C_H$ are the same when $de$ is even and $r \neq 3$ or when $r=3$ and
$e \in \{1,3\}$ or when $r=2$ and $d=e=1$.

For the exceptional types, $G_{25}$ is the only case where the commuting orbits under $N_H$ and $C_H$ are not the same.
The only exceptional types where the non commuting orbits under $N_H$ and $C_H$ are not the same are $G_4$, the second (named after {\sf GAP})
class of hyperplanes of $G_6$, the first (named after {\sf GAP}) class of hyperplanes of $G_{13}$ and the third (named after {\sf GAP}) class of hyperplanes of $G_{15}$
(see {\sf GAP} instructions in Remark~\ref{rem-gap}).
\end{prop-eng}

\begin{demo-eng} From Lemma~\ref{lem-abelianization}, it is enough to show that every linear character of $p^{-1}(C_H)$ with values in $\CC^\times$
extends to $p^{-1}(N_H)$. But the group of linear characters of $p^{-1}(C_H)$ is generated by the $\exp(z \ph_\Cc)$ for $z \in \CC$ and $\Cc$ an orbit
of $\Hscr$ under $C_H$. So it suffices to show that $\ph_\Cc$ extends to $p^{-1}(N_H)$.

Since the orbits of $\Hscr$ under $C_H$ and $N_H$ are the same, for every $\Cc \in \Hscr/C_H$, there exists $n \in \mathbb{N}^*$ such that
${\al_\Cc}^n$ is invariant under $N_H$ (see Lemma~\ref{lem-invariant} for the definition of $\al_\Cc$).
Then Remark~\ref{rem-prolongement} shows that
$$\psi_\Cc : \ga \in p^{-1}(N_H) \longmapsto \frac{1}{n}\int_{{\al_\Cc}^n \circ \ga} \frac{\textrm{d} z}{z} \in \QQ$$
is a well defined linear character of $p^{-1}(N_H)$ extending $\ph_\Cc$.

Proposition~\ref{prop-abelianisation} applied to $C_H$ ensures us that $p^{-1}(C_H)^{\textrm{\scriptsize ab}}$ is the free abelian group
generated by $({s_{\{H\}}^{e_{\{H\}}}})^{\textrm{\scriptsize ab}}$, $(s_\Cc)^{\textrm{\scriptsize ab}}$ for $\Cc \in \Hscr'_H/C_H$
and $(s_{\Cc}^{e_\Cc})^{\textrm{\scriptsize ab}}$ for $\Cc \in (\Hscr\setminus \Hscr_H)/C_H$. Moreover, we have $\widetilde{w}^{f_H} \in p^{-1}(C_H)$
and, thanks to Remark~\ref{rem-prolongement},
$$\ph_{\{H\}}(\widetilde{w}^{f_H}) = \frac{1}{f_H} \rho'(\widetilde{w}^{f_H}) = 1$$
We then deduce that $\widetilde{w}^{f_H}$, $(s_\Cc)^{\textrm{\scriptsize ab}}$ for $\Cc \in \Hscr'_H/C_H$
and $(s_{\Cc}^{e_\Cc})^{\textrm{\scriptsize ab}}$ for $\Cc \in (\Hscr\setminus \Hscr_H)/C_H$ is a basis for $p^{-1}(C_H)^{\textrm{\scriptsize ab}}$.
The short exact sequence
$$\xymatrix{1 \ar[r] & p^{-1}(C_H)^{\textrm{\scriptsize ab}} \ar^{j^{\textrm{\scriptsize ab}}}[r]
    & p^{-1}(N_H)^{\textrm{\scriptsize ab}} \ar^-{(rp)^{{\textrm{\scriptsize ab}}}}[r] & \UU_{f_H} \ar[r]& 1}$$
gives the result.
\end{demo-eng}

\section{An Application to Cohomology}\label{sec-application-cohom}

In this section, we apply the preceding constructions and results to obtain a group cohomology result.
Specifically, the derived subgroup of $P$ is normal in $B$, so we obtain the following short exact sequence
\begin{equation}\label{pbw}
\xymatrix{1\ar[r]  & P/[P,P] \ar^{j}[r] & B/[P,P] \ar^-{p}[r] & W \ar[r] & 1}
\end{equation}
which induces a structure of $W$-module on $\Pab$. By a classical result on hyperplanes arrangements (see~\cite{orlik-solomon} for example),
the $W$-module $\Pab$ is nothing else that the permutation module $\ZZ\Hyp$ and this section describes the extension~(\ref{pbw})
as an element of $H^2(W,\ZZ\Hscr)$ using methods of low-dimensional cohomology.

The rationale breaks down into three steps and each step consists of a translation of a standard isomorphism between cohomology groups in terms of group extensions.

\begin{ri}{(iii)}
\item We decompose $\Hyp$ into orbits under $W$ : $\Hyp = \sqcup\, \Cc$ and uses the isomorphism
$$H^2(W,\ZZ\Hyp) = \bigoplus_{\Cc \in \Hyp/W} H^2(W,\ZZ\Cc)$$
\item In each orbit, we choose a hyperplane $H_\Cc$ and write $\ZZ\Cc = \Ind{N_\Cc}{W} \ZZ$ where $N_\Cc$ is the stabilizer of $H_\Cc$.
Shapiro's isomorphism lemma (see~\cite[Proposition III.6.2]{brown}) then gives us
$$H^2(W,\ZZ\Hyp) = \bigoplus_{\Cc \in \Hyp/W} H^2(N_\Cc,\ZZ)$$
\item The short exact sequence~$\xymatrix{0 \ar[r]& \ZZ \ar[r]& \QQ \ar[r]& \QQ/\ZZ \ar[r]& 0}$ of $N_\Cc$-modules gives us a long exact sequence in
cohomology. Since $H^1(N_{\Cc}, \QQ) = H^2(N_{\Cc}, \QQ) = 0$ (since $|N_{\Cc}|$ is invertible in $\QQ$), we obtain $H^2(N_\Cc,\ZZ) = H^{1}(N_\Cc,\QQ/\ZZ)$ and
$$H^2(W,\ZZ\Hyp) = \bigoplus_{\Cc \in \Hyp/W} H^1(N_{\Cc},\QQ/\ZZ) = \bigoplus_{\Cc \in \Hyp/W} \Hom{\textrm{gr.}}{N_{\Cc}}{\QQ/\ZZ}\,.$$
\end{ri}

The results of this section are the following proposition and corollary.

\begin{prop-eng}[Description]\label{prop-description} Under the isomorphism
$$H^2(W,\ZZ\Hyp) = \bigoplus_{\Cc \in \Hyp/W} \Hom{\textrm{gr.}}{N_{\Cc}}{\QQ/\ZZ}$$
the extension~(\ref{pbw}) correspond to the family of morphisms $(r_\Cc: N_\Cc\ra \QQ/\ZZ)_{\Cc \in \Hyp/W}$ where $r_\Cc$ is the
natural linear character of $N_\Cc$ (we identify $\UU_{f_{H_\Cc}}$ with a subgroup of $\QQ/\ZZ$ via the exponential map).
\end{prop-eng}

The next corollary is a trivial consequence of Proposition~\ref{prop-description} and generalizes a result of Digne~\cite[5.1]{digne} for the case of Coxeter groups. 

\begin{cor-eng}[Order in $H^2(W,\ZZ\Hyp)$]\label{cor-order-cocycle}
In particular, since the order of $r_\Cc$ is $f_{H_\Cc}$, we obtain that the order of the extension~(\ref{pbw}) is
$\kappa(W)=\textrm{lcm}(f_{H_\Cc}, \Cc \in \Hyp/W)$ (this integer $\kappa(W)$ was first introduced in~\cite{marin}).
\end{cor-eng}

The rest of the section is devoted to the proof of Proposition~\ref{prop-description} : one subsection for each of the three steps.

\subsection{First step : splitting into orbits}

The isomorphism
$$H^2(W,\ZZ\Hyp) = \bigoplus_{\Cc \in \Hscr/W} H^2(W,\ZZ\Cc)$$
\noindent is simply given by applying the various projections $p_\Cc : \ZZ\Hscr \ra \ZZ\Cc$ to a $2$-cocycle with values in $\ZZ\Hscr$ where
$$p_\Cc : \sum_{H \in \Hscr} \lambda_H H \longmapsto \sum_{H \in \Cc} \lambda_H H\,.$$

To give a nice expression of the corresponding extensions, we need the following lemma.

\begin{lem-eng}[Extension and direct sum]\label{lem-ext-directsum} Let $G$ be a group, $X= Y \oplus Z$ a direct sum of $G$-modules and
$$\extensionzu{X}{E}{G}{u}{v}$$
\noindent an extension of $G$ by $X$. We denote by $q : X \ra Y$ the first projection and $\ph$ the class of the extension $E$ in $H^2(G,X)$.
The extension associated to $q(\ph)$ is
$$\extensionzu{Y}{E/Z}{G}{}{}$$
\end{lem-eng}

\begin{demo-eng} Let us denote $\theta : E \ra E/Z$ the natural surjection and $i : Y \ra Y \oplus Z$ the natural map.
Let us first remark that $Z$ is normal in $E$ since $Z$ is stable by the action of $G$. Since $v$ is trivial on $Z$, then it induces a
group homomorphism $\widetilde{v} : E/Z \ra G$ whose kernel is $X/Z=Y$. Thus the sequence
\begin{equation}\label{ext-quotient}
\extensionzu{Y}{E/Z}{G}{\theta u i}{\widetilde{v}}
\end{equation}
\noindent is an exact one.

Let us choose $s : G \ra E$ a set-theoretic section of $v$. A $2$-cocycle associated to the extension $E$ is then given by $\ph : G^2 \ra X$
such that $u(\ph(g,g'))= s(g)s(g')s(gg')^{-1}$ for $g,g' \in G$.

Moreover, the kernel of $\theta$ is $Z$, so $\theta u i (q(\ph(g,g'))) = \theta u(\ph(g,g'))= \theta s(g)\theta s(g') \theta(s(gg'))^{-1}$.
Since $\theta s$ is a section for~(\ref{ext-quotient}), a $2$-cocycle associated to~(\ref{ext-quotient}) is given by $(g,g') \mapsto q(\ph(g,g'))$.
\end{demo-eng}

For $\Cc \in \Hyp/W$, we denote $B_\Cc$ the quotient group
$$B_\Cc= B/\langle [P,P], {s_{H,\ga}}^{e_H},\quad H \notin \Cc \rangle$$
\noindent Lemma~\ref{lem-ext-directsum} tells us that the extension~(\ref{pbw}) is equivalent to the family of extensions
\begin{equation}\label{pbw-loc}
\extensionzu{\ZZ\Cc}{B_\Cc}{W}{j_\Cc}{p_\Cc}
\end{equation}
\noindent for $\Cc \in \Hyp/W$.

\subsection{Second step : the induction argument}

\vskip1ex
In each orbit $\Cc \in \Hyp/W$, we choose a hyperplane $H_\Cc \in \Cc$ and write $\ZZ\Cc = \Ind{N_{\Cc}}{W} \ZZ$ where $N_{\Cc}\subset W$ is the stabilizer of $H_\Cc$.
The Shapiro's isomorphism lemma~\cite[Proposition III.6.2]{brown} shows that $H^2(W,\ZZ\Cc) = H^2(N_\Cc,\ZZ)$.
Exercise~III.8.2 of~\cite{brown} tells us that in term of $2$-cocycles the Shapiro's isomorphism is described as follow
$$S: (\ph: G^2 \rightarrow \ZZ\Cc) \longmapsto (f_\Cc \circ \ph : {N_\Cc}^2 \rightarrow \ZZ)$$
\noindent where $f_\Cc : \ZZ\Cc \ra \ZZ$ is the projection onto the $H_\Cc$-component.

Decomposing the Shapiro's isomorphism into the following two steps
$$(\ph: G^2 \rightarrow \ZZ\Cc)\longmapsto (\ph: {N_\Cc}^2 \rightarrow \ZZ\Cc) \longmapsto (f_\Cc \circ \ph : {N_\Cc}^2 \rightarrow \ZZ)\,,$$
allows us to interpret it in terms of group extensions. Exercice~IV.3.1.(a) of~\cite{brown} gives a description of the first step :
the corresponding extension is given by
$$\extensionzu{\ZZ\Cc}{{p_\Cc}^{-1}(N_\Cc)}{N_\Cc}{}{}$$
since ${p_\Cc}^{-1}(N_\Cc)$ is the fiber product of $B_\Cc$ and $N_\Cc$ over $W$. Moreover, since $f_\Cc$ is a split surjection as a $N_\Cc$-module map, Lemma~\ref{lem-ext-directsum} gives us the following extension
$$\extensionzu{\ZZ}{{p_\Cc}^{-1}(N_\Cc)/\langle {s_{H,\ga}}^{e_H},\  H \in \Cc \setminus \{ H_\Cc \} \rangle}{N_\Cc}{}{}$$

Finally, the extension~(\ref{pbw}) is equivalent to the family of extensions
\begin{equation}\label{ext-induit}
\extensionzu{\ZZ}{B'_\Cc}{N_\Cc}{}{p_\Cc}
\end{equation}
\noindent where $B'_\Cc= p^{-1}(N_\Cc)/\langle [P,P], {s_{H,\ga}}^{e_H}, \ H \neq H_\Cc \rangle$ and $\Cc \in \Hscr/W$.

\subsection{Third step : linear character}

\vskip1ex
For the third step, we use results and notations of Section~\ref{sec-reflection-group-case} and Section~\ref{sec-braid-group-case}.
Let us consider the group homomorphism $\rho_\Cc : p^{-1}(N_\Cc) \ra \QQ$ of Definition~\ref{dfn-morphism}.
Since it is trivial on $\langle [P,P], {s_{H,\ga}}^{e_H}, \ H \neq H_\Cc \rangle$,
it induces a group homomorphism from $B'_\Cc$ to $\QQ$ still denoted $\rho_\Cc$. Moreover, since $\rho_\Cc({s_{H_\Cc,\ga}}^{e_{H_\Cc}}) = 1$,
Proposition~\ref{prop-lin-caractere} gives the following commutative diagram
$$\xymatrix{0 \ar[r] & \ZZ \ar[r] \ar@{=}[d]& \QQ \ar[r] & \QQ/\ZZ \ar[r]  & 0\\
        0 \ar[r] & \ZZ \ar[r] & B'_\Cc \ar^{p_\Cc}[r] \ar_{\rho_\Cc}[u] & N_\Cc \ar[r] \ar_{r_\Cc}[u] & 1}$$
\noindent Exercises~IV.3.2 and~IV.3.3 of~\cite{brown} tell us precisely that the group homomorphism corresponding to~(\ref{ext-induit}) is $r_\Cc$.
So the extension~(\ref{pbw}) is equivalent to the family $(r_\Cc)_{\Cc \in \Hscr/W}$. This concludes the proof of Proposition~\ref{prop-description}.

\section{Tables}\label{sec-table}

\subsection{The infinite series}

In this subsection, we give tables for the orbits of the hyperplanes of $G(de,e,r)$ under the centralizer of a reflection
and under the parabolic subgroup associated to the line of the reflection. We also gives tables for the values of $f_H$ and the index of ramification. 
So let us consider the complex reflection group $G(de,e,r)$ acting on $\CC^r$ with canonical basis $(e_1,\ldots, e_r)$. 
The standard point of $\CC^r$ is denoted by $(z_1,\ldots, z_r)$.

The hyperplanes of $G(de,e,r)$ are precisely the $H_i = \{z_i=0\}$ for $i \in \{1,\ldots, r\}$ (when $d > 1$) and the
$H_{i,j,\zeta} = \{z_i=\zeta z_j\}$ for $i<j$ and $\zeta \in \UU_{de}$ (when $r \geq 2$). They split in general into two conjugacy classes under
$G(de,e,r)$ whose representant may be chosen as follow $H_1$ and $H_{1,2,1}$.

Let us continue with more notations. For every triple of integers $d,e,r$, we denote by $\pi : G(de,e,r) \ra \UU$ the following group morphism :
for $g \in G(de,e,r)$, $\pi(g)$ is the product of the nonzero coefficients of the monomial matrix $g$. When $e$ is even, we denote by $e' = e/2$.
We denote by $e''=e/\gcd(e,3)$ and by $P$ the set of elements of $\UU_{de}$ with strictly positive imaginary part.

\subsubsection{The case of the hyperplane $H_1 = \{z_1=0\}$}

We then have $d > 1$. The stabilizer $N$ of $H_1$ is described by
$$N=\{(\al,g), \quad g \in G(de,1,r-1),\ \al \in \UU_{de}, \ (\pi(g)\al)^d = 1\}$$
and the pointwise stabilizer $C$ of $D_1={H_1}^{\perp} = \CC e_1$ is described by $C=G(de,e,r-1)$

The following table gives the orbits of the hyperplanes under $N$ and $C$.

$$\begin{array}{|c|c|c|c|}\hline  &&N  & C \\ \hline
r \geq 4 & \textrm{commuting} & H_1 & H_1\\
&& \{H_i,\ i\geq 2\} & \{H_i,\ i\geq 2\}\\
&& \{ H_{i,j,\zeta}, \ i\neq j \geq 2,\ \zeta \in \UU_{de} \} & \{ H_{i,j,\zeta}, \ i\neq j \geq 2,\ \zeta \in \UU_{de} \}\\
& \textrm{non commuting} & \{H_{1,j,\zeta},\ j \geq 2,\ \zeta \in \UU_{de}\} & \{H_{1,j,\zeta},\ j \geq 2,\ \zeta \in \UU_{de}\}\\ \hline
r=3 & \textrm{commuting} & H_1 & H_1 \\
\textrm{and} && \{ H_2,H_3 \} & \{ H_2,H_3 \}\\
e \textrm{ odd} && \{ H_{2,3,\zeta},\ \zeta \in \UU_{de} \}& \{ H_{2,3,\zeta},\ \zeta \in \UU_{de} \}\\
& \textrm{non commuting} & \{ H_{1,i,\zeta},\ i=2,3,\ \zeta \in \UU_{de}\} & \{ H_{1,i,\zeta},\ i=2,3,\ \zeta \in \UU_{de}\}\\ \hline
r=3 & \textrm{commuting} & H_1 & H_1 \\
\textrm{and} && \{ H_2,H_3 \} & \{ H_2,H_3 \}\\
e \textrm{ even} && \{ H_{2,3,\zeta},\ \zeta \in \UU_{de} \}& \{ H_{2,3,\zeta},\ \zeta \in c \} \textrm{ for } c \in \UU_{de}/\UU_{de'}\\
& \textrm{non commuting} & \{ H_{1,i,\zeta},\ i=2,3,\ \zeta \in \UU_{de}\} & \{ H_{1,i,\zeta},\ i=2,3,\ \zeta \in \UU_{de}\}\\ \hline
r=2 & \textrm{commuting} & H_1& H_1 \\
\textrm{and} && H_2 & H_2 \\
e \textrm{ odd} &\textrm{non commuting}& \{ H_{1,2,\zeta},\ \zeta \in \UU_{de} \} &  \{ H_{1,2,\zeta},\ \zeta \in c \} \textrm{ for } c \in \UU_{de}/\UU_d\\ \hline
r=2 & \textrm{commuting} & H_1& H_1 \\
\textrm{and} && H_2 & H_2 \\
e \textrm{ even} &\textrm{non commuting}& \{ H_{1,2,\zeta},\ \zeta \in c \} \textrm{ for } c \in \UU_{de}/\UU_{de'} &
\{ H_{1,2,\zeta},\ \zeta \in c \} \textrm{ for } c \in \UU_{de}/\UU_d\\ \hline
r=1 &\textrm{commuting} & H_1 & H_1\\ \hline
\end{array}$$





\subsubsection{The hyperplane $H_{1,2, \exp(2i\pi/de)}$ with $r=2$}

Let us define $\zeta = \exp(2i\pi/de)$. The reflection of $G(de,e,2)$ with hyperplane $H_{1,2,\zeta}$ is
$$s = \left[\begin{array}{@{}c@{\ }c@{}} & \zeta \\ \zeta^{-1} & \end{array} \right]$$
The line of $s$ is $D=\CC(e_2-\zeta e_1)$. The centralizer of $s$ is given by
$$N= \left\{ d_\lambda=\left[\begin{array}{@{}c@{\ }c@{}} \lambda &  \\  &\lambda \end{array} \right],
        t_\lambda=\left[\begin{array}{@{}c@{\ }c@{}} & \lambda \zeta \\ \lambda \zeta^{-1} & \end{array} \right],\  \lambda^{2d}=1 \right\}\,.$$
The eigenvalue of $t_\lambda$ on $D$ is $\lambda$ whereas the eigenvalue of $d_\lambda$ on $D$ is $-\lambda$.
So the parabolic subgroup $C$ associated to $D$ is $C= \left\{ \id, -s \right\}$.
The orbits of hyperplane under $C$ and $N$ are the same. The commuting ones are $\{H_{1,2,\zeta}\}$ and $\{H_{1,2,-\zeta}\}$.
The non commuting ones are $\{H_1,H_2\}$ and $\{H_{1,2,\mu}, H_{1,2,\zeta^2\mu^{-1}}\}$ for $\mu \in \UU_{de} \setminus \{ \pm \zeta \}$.

\subsubsection{The hyperplane $H_{1,2,1} = \{ z_1=z_2\}$} 

We then have $r \geq 2$. The reflection of $G(de,e,r)$ with hyperplane $H_{1,2,1}$ is the transposition $\tau_{12}$ swapping $1$ and $2$.
Since the elements of $G(de,e,r)$ are monomial matrices, an element of $G(de,e,r)$ commuting with $\tau_{12}$ stabilizes the subspace
spanned by $e_1$ and $e_2$. Thus the stabilizer $N$ of $H_{1,2,1}$ is given by
$$N = \left\{d_{\lambda,g}=\left[\begin{array}{@{}c@{\ }c@{\ }|c@{}} \lambda && \\[-0.5ex]  &\lambda&  \\[-0.7ex]\hline  &&g \end{array}\right],
 t_{\lambda,g}=\left[\begin{array}{@{}c@{\ }c@{\ }|c@{}}  &\lambda& \\[-0.5ex]  \lambda&&  \\[-0.7ex]\hline  &&g \end{array}\right], \quad
        \lambda \in \UU_{de},\ g \in G(de,1,r-2),\ (\pi(g)\lambda^2)^d = 1 \right\}$$
The line of $\tau_{12}$ is $\CC(e_1-e_2)$. So the eigenvalue of $d_{\lambda,g}$ on $\CC(e_1-e_2)$ is $\lambda$ whereas the eigenvalue of $t_{\lambda,g}$
on $\CC(e_1-e_2)$ is $-\lambda$. Thus, when $de$ is odd, the parabolic subgroup associated to the line $\CC(e_1 - e_2)$ is given by
$$C = \left\{ \left[\begin{array}{@{}c@{\ }c@{\ }|c@{}} 1 && \\[-0.5ex]  &1&  \\[-0.7ex]\hline  &&g \end{array}\right],\quad g \in G(de,e,r-2) \right\}$$
and when $de$ is even, the parabolic subgroup associated to the line $\CC(e_1 - e_2)$ is given by
$$C = \left\{ \left[\begin{array}{@{}c@{\ }c@{\ }|c@{}} 1 && \\[-0.5ex]  &1&  \\[-0.7ex]\hline  &&g \end{array}\right],
\left[\begin{array}{@{}c@{\ }c@{\ }|c@{}}  &-1& \\[-0.5ex]  -1&&  \\[-0.7ex]\hline  &&g \end{array}\right], \quad g \in G(de,e,r-2) \right\}$$

The following table gives the orbits the hyperplanes of $G(de,e,r)$ under $N$ and $C$.

$$\begin{array}{|c|c|c|c|}\hline  &&N  & C \\ \hline
r\geq 5 &\textrm{commuting}& H_{1,2,1}& H_{1,2,1}\\
\textrm{and} && \{ H_i,\ i\geq 3 \} \textrm{ if } d > 1 & \{ H_i,\ i\geq 3 \} \textrm{ if } d > 1 \\
de \textrm{ odd} && \{ H_{i,j,\zeta},\  i\neq j \geq 3,\ \zeta \in \UU_{de}\}& \{ H_{i,j,\zeta},\  i\neq j \geq 3,\ \zeta \in \UU_{de}\}\\
&\textrm{non commuting}& \{H_1,H_2\} \textrm{ if } d > 1& \{H_1\}, \{H_2\} \textrm{ if } d > 1\\
&& \{H_{1,2,\zeta}, H_{1,2,\zeta^{-1}}\} \textrm{ for } \zeta \in P   & \{H_{1,2,\zeta}\} \textrm{ for } \zeta \in \UU_{de} \\
&& \{H_{1,j,\zeta}, H_{2,j,\zeta},\ j \geq 3,\ \zeta \in \UU_{de} \} & \{H_{1,j,\zeta}, \ j \geq 3,\ \zeta \in \UU_{de} \}\\
&&&    \{H_{2,j,\zeta}, \ j \geq 3,\ \zeta \in \UU_{de} \}\\ \hline
r\geq 5 &\textrm{commuting}& H_{1,2,1}& H_{1,2,1}\\
\textrm{and}& & H_{1,2,-1}  & H_{1,2,-1} \\
de \textrm{ even}&& \{ H_i,\ i\geq 3 \} \textrm{ if } d > 1 & \{ H_i,\ i\geq 3 \} \textrm{ if } d > 1 \\
&& \{ H_{i,j,\zeta},\  i\neq j \geq 3,\ \zeta \in \UU_{de}\}& \{ H_{i,j,\zeta},\  i\neq j \geq 3,\ \zeta \in \UU_{de}\}\\
&\textrm{non commuting}& \{H_1,H_2\} \textrm{ if } d > 1& \{H_1,H_2\} \textrm{ if } d > 1\\
&& \{H_{1,2,\zeta}, H_{1,2,\zeta^{-1}}\} \textrm{ for } \zeta \in P   & \{H_{1,2,\zeta}, H_{1,2,\zeta^{-1}}\} \textrm{ for } \zeta \in P \\
&& \{H_{1,j,\zeta}, H_{2,j,\zeta},\ j \geq 3,\ \zeta \in \UU_{de} \} & \{H_{1,j,\zeta}, H_{2,j,\zeta},\ j \geq 3,\ \zeta \in \UU_{de} \}\\ \hline
r=4 &\textrm{commuting}&H_{1,2,1}& H_{1,2,1}\\
\textrm{and}  && \{ H_i,\ i\geq 3 \} \textrm{ if } d > 1 & \{ H_i,\ i\geq 3 \} \textrm{ if } d > 1 \\
d,e \textrm{ odd}&& \{ H_{i,j,\zeta},\  i\neq j \geq 3,\ \zeta \in \UU_{de}\}& \{ H_{i,j,\zeta},\  i\neq j \geq 3,\ \zeta \in \UU_{de}\}\\
& \textrm{non commuting} & \{ H_1,H_2\} \textrm{ if } d > 1& \{H_1 \}, \{H_2\} \textrm{ if } d > 1 \\
&& \{H_{1,2,\zeta}, H_{1,2,\zeta^{-1}}\} \textrm{ for } \zeta \in P   & \{H_{1,2,\zeta}\} \textrm{ for } \zeta \in \UU_{de} \\
&& \{H_{1,j,\zeta}, H_{2,j,\zeta},\ j \geq 3,\ \zeta \in \UU_{de} \} & \{H_{1,j,\zeta}, \ j \geq 3,\ \zeta \in \UU_{de} \}\\
&&&    \{H_{2,j,\zeta}, \ j \geq 3,\ \zeta \in \UU_{de} \}\\ \hline
r=4 &\textrm{commuting}&H_{1,2,1}& H_{1,2,1}\\
\textrm{and} && H_{1,2,-1} & H_{1,2,-1} \\
e \textrm{ odd} && \{ H_i,\ i\geq 3 \} \textrm{ if } d > 1 & \{ H_i,\ i\geq 3 \} \textrm{ if } d > 1 \\
d \textrm{ even}&& \{ H_{i,j,\zeta},\  i\neq j \geq 3,\ \zeta \in \UU_{de}\}& \{ H_{i,j,\zeta},\  i\neq j \geq 3,\ \zeta \in \UU_{de}\}\\
&\textrm{non commuting}& \{H_1,H_2\} \textrm{ if } d > 1& \{H_1,H_2\} \textrm{ if } d > 1\\
&& \{H_{1,2,\zeta}, H_{1,2,\zeta^{-1}}\} \textrm{ for } \zeta \in P   & \{H_{1,2,\zeta}, H_{1,2,\zeta^{-1}}\} \textrm{ for } \zeta \in P \\
&& \{H_{1,j,\zeta}, H_{2,j,\zeta},\ j \geq 3,\ \zeta \in \UU_{de} \} & \{H_{1,j,\zeta}, H_{2,j,\zeta},\ j \geq 3,\ \zeta \in \UU_{de} \}\\ \hline
r=4 & \textrm{commuting} &H_{1,2,1}& H_{1,2,1}\\
\textrm{and}& & H_{1,2,-1}  & H_{1,2,-1} \\
e \textrm{ even}&& \{ H_i,\ i\geq 3 \} \textrm{ if } d > 1 & \{ H_i,\ i\geq 3 \} \textrm{ if } d > 1 \\
&& \{ H_{i,j,\zeta},\  i\neq j \geq 3,\ \zeta \in \UU_{de'}\}& \{ H_{i,j,\zeta},\  i\neq j \geq 3,\ \zeta \in \UU_{de'}\}\\
&& \{ H_{i,j,\zeta},\  i\neq j \geq 3,\ \zeta \in \UU_{de} \setminus \UU_{de'}\}& \{ H_{i,j,\zeta},\  i\neq j \geq 3,\ \zeta \in \UU_{de} \setminus  \UU_{de'}\}\\
&\textrm{non commuting}& \{H_1,H_2\} \textrm{ if } d > 1& \{H_1,H_2\} \textrm{ if } d > 1\\
&& \{H_{1,2,\zeta}, H_{1,2,\zeta^{-1}}\} \textrm{ for } \zeta \in P   & \{H_{1,2,\zeta}, H_{1,2,\zeta^{-1}}\} \textrm{ for } \zeta \in P \\
&& \{H_{1,j,\zeta}, H_{2,j,\zeta},\ j \geq 3,\ \zeta \in \UU_{de} \} & \{H_{1,j,\zeta}, H_{2,j,\zeta},\ j \geq 3,\ \zeta \in \UU_{de} \}\\ \hline
r=3 &\textrm{commuting}& H_{1,2,1}& H_{1,2,1}\\
\textrm{and} && H_3 \textrm{ if } d > 1 & H_3 \textrm{ if } d > 1 \\
de \textrm{ odd} & \textrm{non commuting} & \{ H_1,H_2 \} \textrm{ if } d > 1& \{ H_1\}, \{H_2 \} \textrm{ if } d > 1\\
&& \{H_{1,2,\zeta}, H_{1,2,\zeta^{-1}}\} \textrm{ for } \zeta \in P   & \{H_{1,2,\zeta}\} \textrm{ for } \zeta \in \UU_{de} \\
&&\{H_{1,3,\zeta}, H_{2,3,\zeta},\ \zeta \in c \}  \textrm{ for } c \in \UU_{de}/\UU_{de''}&
U_c = \{H_{1,3,\zeta}, \, \zeta \in c \} \textrm{ for } c \in \UU_{de}/\UU_d \\
&&& U'_c = \{H_{2,3,\zeta}, \, \zeta \in c \} \textrm{ for } c \in \UU_{de}/\UU_d\\ \hline
r=3 &\textrm{commuting}& H_{1,2,1}& H_{1,2,1}\\
\textrm{and} && H_{1,2,-1}& H_{1,2,-1}\\
de \textrm{ even} && H_3 \textrm{ if } d > 1 & H_3 \textrm{ if } d > 1 \\
& \textrm{non commuting} & \{ H_1,H_2 \} \textrm{ if } d > 1& \{ H_1,H_2 \} \textrm{ if } d > 1\\
&& \{H_{1,2,\zeta}, H_{1,2,\zeta^{-1}}\} \textrm{ for } \zeta \in P   & \{H_{1,2,\zeta}, H_{1,2,\zeta^{-1}}\} \textrm{ for } \zeta \in P  \\
&&\{H_{1,3,\zeta}, H_{2,3,\zeta},\, \zeta \in c \}  \textrm{ for } c \in \UU_{de}/\UU_{de''}&
\{H_{1,3,\zeta}, H_{2,3,-\zeta}, \, \zeta \in c \} \textrm{ for } c \in \UU_{de}/\UU_d\\ \hline
r=2 &\textrm{commuting}& H_{1,2,1}& H_{1,2,1}\\
\textrm{and} & \textrm{non commuting}& \{ H_1,H_2\} \textrm{ if } d > 1& \{H_1\}, \{H_2\} \textrm{ if } d > 1\\
de \textrm{ odd } && \{H_{1,2,\zeta}, H_{1,2,\zeta^{-1}}\} \textrm{ for } \zeta \in P   & \{H_{1,2,\zeta}\} \textrm{ for } \zeta \in \UU_{de} \\ \hline
r=2 &\textrm{commuting}& H_{1,2,1}& H_{1,2,1}\\
\textrm{and}&& H_{1,2,-1}& H_{1,2,-1}\\
de \textrm{ even } & \textrm{non commuting}& \{ H_1,H_2\} \textrm{ if } d > 1& \{H_1, H_2\} \textrm{ if } d > 1\\
&& \{H_{1,2,\zeta}, H_{1,2,\zeta^{-1}}\} \textrm{ for } \zeta \in P   & \{H_{1,2,\zeta}, H_{1,2,\zeta^{-1}}\} \textrm{ for } \zeta \in P \\ \hline
\end{array}$$

\subsubsection{Value for $f_H$ and the index of ramification}
The computations of the preceding paragraphs gives the following table. We set $\zeta = \exp(2i\pi/de)$.

$$\begin{array}{|r@{\quad}c@{\quad}c|c|c|c|} \hline &H&&f_H& e_H& f_H/e_H \\  \hline
x_1=0 & r=1  &&  d&d &1 \\ \hline
x_1=0 & r\geq 2 &d \neq 1& de& d & e \\ \hline
x_2 =x_1  &r\geq 3& de  \textrm{ odd}  &2de&2& de\\ \hline
x_2=x_1 &r\geq 3& de  \textrm{ even} & de& 2& de/2\\  \hline
x_2=x_1 &r=2& e \textrm{ odd \  and \ } d \textrm{ even}&d&2&d/2\\ \hline
x_2=x_1 &r=2& e \textrm{ even\  or \ } d \textrm{ odd}&2d&2&d\\ \hline
x_2=\zeta x_1 &r=2& e \textrm{ even} &2d&2&d\\ \hline
\end{array}$$
We obtain the following errata for the proposition 6.1 of~\cite{marin}. Let us consider $r \geq 2$. 
For $W=G(de,e,r)$, we have $\kappa(W) = 2de$ if $de$ is odd and $r \geq 3$.
We have $\kappa(W) = de$ if ($d \neq 1$ and $r=2$) or ($r \geq 3$ and $de$ even).

\subsection{Exceptional types}

With the package {\sf CHEVIE} of {\sf GAP} \cite{gap}\cite{chevie} (see instructions in Remark~\ref{rem-gap}),
we obtain the following table for the values of $e_H,f_H$ and $f_H/e_H$ for the hyperplanes of the exceptional reflection groups.
In particular, the only non Coxeter groups with only unramified hyperplanes are $G_8,G_{12}$ and $G_{24}$.
This table can also easily be obtained from the table of~\cite{marin} for the value of $\kappa(W)$.
The first, fifth and ninth columns stand for the number of the group in the Shephard and Todd classification.

\begin{table}[!h]\label{table-fh}
$\begin{array}{|c|c|c|c||c|c|c|c||c|c|c|c||} \hline  ST &e_H& f_H & f_H/e_H &ST &e_H& f_H & f_H/e_H &ST &e_H& f_H & f_H/e_H\\ \hline
4 &3& 6& 2&%
16&5&10&2&%
28&2,2&2,2&1,1\\
5 &3,3& 6,6 & 2,2&%
17&2,5&20,20&10,4&%
29&2&4&2\\
6 &2,3&4,12 & 2,4&%
18&3,5&30,30&10,6&%
30&2&2&1\\
7 &2,3,3&12,12,12 & 6,4,4&%
19&2,3,5&60,60,60&30,20,12&%
31&2&4&2\\
8 &4&4&1&%
20&3&6&2&%
32&3&6&2\\
9 &2,4 & 8,8&4,2&%
21&2,3 & 12,12& 6,4&%
33&2&6&3\\
10&3,4& 12,12&4,3&%
22&2&4&2&%
34&2&6&3\\
11&2,3,4&24,24,24&12,8,6&%
23&2&2&1&%
35&2&2&1\\
12&2&2&1&%
24&2&2&1&%
36&2&2&1\\
13&2,2&8,4&4,2&%
25&3&6&2&%
37&2&2&1\\
14&2,3&6,6&3,2&%
26&3,2&6,6&2,3\\
15&2,3,2&12,12,24&6,4,12&%
27&2&6&3\\
\end{array}$
\caption{Values for the ramification index}
\end{table}

\begin{rem-eng}[{\bf\sf GAP} Instructions]\label{rem-gap}
\begin{verbatim}
RequirePackage("chevie");
for j in [4..37]
do
	G:=ComplexReflectionGroup(j);
	Print("G"); Print(j); Print("\n");
	Ref:=Set(Reflections(G));
	L:=HyperplaneOrbits(G); T:=Length(L);
	Print("There are "); Print(T); Print(" classes of hyperplanes\n");
	for k in [1..T]
	do

#	Construction of the subgroup N_H and its abelianization which are denoted by N and Nab
		RefIndice:=L[k].s;
		s:=Reflections(G)[RefIndice];
		H:=Group(s);
        e:=Order(G,s);
		N:=Centralizer(G,s);
		Nab:=CommutatorFactorGroup(N);
		nab:=Size(Nab);

#	Construction of the subgroup C_H and its abelianization which are denoted by C and Cab
		r:=[];
		for i in [1..Size(Ref)]
		do
			if Comm(Ref[i],s) = G.identity and (not IsSubset(H,[Ref[i]]))
				then r:=Concatenation(r,[Ref[i]]);
			fi;
		od;
		C:=Subgroup(G,r); f:=Size(FactorGroup(N,C));
		Cab:=CommutatorFactorGroup(C); cab:=Size(Cab);

#	Determination of the orbites of hyperplanes under N and C which commute with s
		tC:=[]; RefOrbComC:=[Size(H)]; tN:=[]; RefOrbComN:=[Size(H)];
		for i in [1..Size(r)]
		do
			if (not IsSubset(tC,[r[i]])) then tC:=Union(tC,ConjugacyClass(C,r[i]));
							Add(RefOrbComC,Size(ConjugacyClass(C,r[i])));
			fi;
			if (not IsSubset(tN,[r[i]])) then tN:=Union(tN,ConjugacyClass(N,r[i]));
							Add(RefOrbComN,Size(ConjugacyClass(N,r[i])));
			fi;
		od;
		
#	Determination of the orbites of hyperplanes under N and C which do not commute with s
		tnC:=[]; RefOrbComnC:=[]; tnN:=[]; RefOrbComnN:=[];
		for i in [1..Size(Ref)]
		do
			if (not Comm(Ref[i],s)=G.identity) and (not IsSubset(tnC,[Ref[i]])) then
tnC:=Union(tnC,ConjugacyClass(C,Ref[i]));
Add(RefOrbComnC,Size(ConjugacyClass(C,Ref[i])));
			fi;
			if (not Comm(Ref[i],s)=G.identity) and (not IsSubset(tnN,[Ref[i]])) then
tnN:=Union(tnN,ConjugacyClass(N,Ref[i])); Add(RefOrbComnN,Size(ConjugacyClass(N,Ref[i])));
			fi;
		od;

# 	The results

Print("For the "); Print(k); Print(" th classe of hyperplanes\n");	
	
	Print("e="); Print(e); Print(" f="); Print(f); Print(" d="); Print(f/e); Print("\n");

	if RefOrbComC = RefOrbComN then
    Print("   There are "); Print(Length(RefOrbComC));
    Print(" classes which commute with s under C and N\n");
    Print("   The commuting orbits under N and C are the same \n");
		else Print("   There are "); Print(Length(RefOrbComC));
    Print(" classes which commute with s under C and ");
    Print(Length(RefOrbComN)); Print(" under N\n");
    Print("   The commuting orbits under N and C are not the same \n");
	fi;

	if RefOrbComnC = RefOrbComnN then
	Print("   There are "); Print(Length(RefOrbComnC));
    Print(" classes which do not commute with s under C and N\n");
    Print("   The non commuting orbits under N and C are the same \n");
else Print("   There are "); Print(Length(RefOrbComnC));
    Print(" classes which do not commute with s under C and ");
    Print(Length(RefOrbComnN)); Print(" under N\n");
    Print("   The non commuting orbits under N and C are not the same \n");
	fi;

	if f*cab-nab = 0 then
    Print("   The abelianized sequence is exact\n");
else Print("   The abelianized sequence is not exact\n");
	fi;
	od;
Print("\n");
od;
\end{verbatim}
\end{rem-eng}

\end{document}